\documentclass[12pt]{amsart}

\usepackage{amssymb,amsmath}
\usepackage{latexsym}

\newtheorem{example}{Example}[section]

\newtheorem{theorem}[example]{Theorem}

\newtheorem{corollary}[example]{Corollary}

\newtheorem{definition}[example]{Definition}
\newtheorem{proposition}[example]{Proposition}

\def\Proof{\noindent \it Proof -- \rm}                                                 
\def\qed{\hspace{3.5mm} \hfill \vbox{\hrule height 3pt depth 2 pt width 2mm}
\bigskip}

\def\sinv{{\rm sinv}}
\def\q{{\bf q}}
\def\t{{\bf t}}
\def\AA{{\mathcal A}}
\def\Sym{{\bf Sym}}

\def\NCSF{{\bf Sym}}

\def\WQSym{{\bf WQSym}}
\def\maj{{\rm maj}}

\def\WC{{\rm WC}}
\def\fatter{\leq}
\def\DesC{{\rm DC}}

\def\<{\langle}
\def\>{\rangle}

\def\RR{{\mathcal R}}
\def\SS{{\mathcal S}}
\def\JJ{{\mathcal J}}
\def\PP{{\mathcal P}}
\def\QQ{{\mathcal Q}}

\def\diag{{\rm diag\,}}

\def\t{{\bf t}}
\def\A{{\sf A}}

\def\GG{{\rm G}}

\def\Des{{\rm Des\,}}

\def\KK{{\mathcal K}}

\def\y{{\bf y}}
\def\x{{\bf x}}

\def\w{{\bf w}}

\def\qbin#1#2{\begin{bmatrix}#1\\ #2 \end{bmatrix}_q}

\newdimen\squaresize
\newdimen\thickness        
                                                    
\def\square#1{\hbox{\vrule width \thickness
   \vbox to \squaresize{\hrule height \thickness\vss                            
      \hbox to \squaresize{\hss#1\hss}
   \vss\hrule height\thickness} 
\unskip\vrule width \thickness} 
\kern-\thickness}                                                            
                               
\def\vsquare#1{\vbox{\square{$\casestyle#1$}}\kern-\thickness}
\def\blank{\omit\hskip\squaresize}

\def\young#1{\vcenter{%
    \squaresize=18pt\thickness=0.5pt\let\casestyle=\relax
    \vbox{\smallskip\offinterlineskip
      \halign{&\vsquare{##}\cr #1}}}}

\def\moyyoung#1{\vcenter{%
    \squaresize=28pt\thickness=0.6pt\let\casestyle=\relax
    \vbox{\smallskip\offinterlineskip
      \halign{&\vsquare{##}\cr #1}}}}

\def\bigyoung#1{\vcenter{%
    \squaresize=38pt\thickness=0.8pt\let\casestyle=\displaystyle
    \vbox{\smallskip\offinterlineskip
      \halign{&\vsquare{##}\cr #1}}}}

\def\boxit#1#2{\setbox1=\hbox{\kern#1{#2}\kern#1}%
\dimen1=\ht1 \advance\dimen1 by #1 \dimen2=\dp1 \advance\dimen2 by #1
\setbox1=\hbox{\vrule height\dimen1 depth\dimen2\box1\vrule}%
\setbox1=\vbox{\hrule\box1\hrule}%
\advance\dimen1 by .4pt \ht1=\dimen1
\advance\dimen2 by .4pt \dp1=\dimen2 \box1\relax}

\def\tbin#1#2{{\begin{bmatrix}#1 \\ #2\end{bmatrix}_t}}
\begin{document}
\title[Noncommutative Hall-Littlewood and Macdonald functions]{On some noncommutative symmetric
functions analogous to Hall-Littlewood and Macdonald polynomials}

\author[J.-C. Novelli, L. Tevlin,
J.-Y. Thibon]%
{Jean-Christophe Novelli, Lenny Tevlin,
and Jean-Yves Thibon}

\address[J.-C Novelli, J.-Y. Thibon] {
Universit\'e Paris-Est Marne-la-Vall\'ee \\
Laboratoire d'Informatique Gaspard-Monge\\
5 Boulevard Descartes \\Champs-sur-Marne \\77454 Marne-la-Vall\'ee cedex 2 \\
France}
\address[L. Tevlin]{
Liberal Studies,\\ New York University\\
726 Broadway\\ New York, N.Y. 10003, U.S.A.}
\thanks{L. Tevlin was partially supported by a grant from the NYU Research Challenge Fund Program.}, 

\dedicatory{To Christophe Reutenauer, on the occasion of his 60th birthday}

\date{\today}



\begin{abstract}
We investigate the connections between various noncommutative
analogues of Hall-Littlewood and Macdonald polynomials,
and
define some new families of noncommutative symmetric functions 
depending on two sequences of parameters.
\end{abstract}

\maketitle

\section{Introduction}

There have been several attempts to define analogues of Hall-Litttlewood and Macdonald
polynomials in
the algebras of noncommutative symmetric functions
($\NCSF$) and of quasi-symmetric functions ($QSym$).
The first analogues of Hall-Littlewood functions were defined by Hivert \cite{Hiv},
who replaced symmetrization by a new operation of quasi-symmetrization in  Littlewood's 
original definition. 
This provided first an interpolation between
the monomial and fundamental bases of quasisymmetric functions
(which are analogues of the monomial and Schur bases of symmetric
functions), by means of a quasi-symmetrizing action of the Hecke
algebra. The construction was then further explored  on the dual side by combinatorial
methods.

It was then shown by Hivert, Lascoux and the third author \cite{HLT}
that these functions admitted a simple and direct combinatorial definition, in which a second
parameter could be introduced so as to give analogues of Macdonald polynomials.
It was also observed that $q$ and $t$ could be replaced by sequences of indeterminates
$q_i$ and $t_i$ in such a way that $q_i=q^i$ and $t_i=t^i$ give back the original
version.

Almost simultaneously, similar but different analogues were
defined by Bergeron and Zabrocki \cite{BZ}.
Their approach was to obtain Macdonald-like
functions from an analogue of the Nabla operator, of which
the original Macdonald functions are the eigenvectors.

More recently \cite{LNT}, it has been shown that many more parameters could be introduced
in the definition of such bases. Actually, one can have a pair of $n\times n$
matrices $(Q_n,T_n)$ for each degree $n$. The main properties established in
\cite{BZ} and \cite{HLT} remain true in this general context, and one recovers
the polynomials introduced in these two papers for appropriate specializations of the matrices.

In the meantime, very different analogues of Hall-Littlewood functions
had been defined by Novelli-Thibon-Williams \cite{NTW}.  
These analogues, which
were based on the  monomial functions
introduced by the second author \cite{Tev07}, had the interesting property
that contrary to  the other versions, the elements of the $(q,t)$-Kostka matrices
were not monomials, but non-trivial polynomials with an interesting
combinatorial interpretation. 
The approach of \cite{NTW} was to build $q$-analogues of product
of complete functions (similar to the $Q'$-version of Hall-Littlewood
functions), such that their expansions on some simple $q$-analogue
of the noncommutative fundamental basis of \cite{Tev07} provided
combinatorial information on permutation tableaux.

Finally, the Hall-Littlewood basis of \cite{Tev} 
was apparently of a different nature. 
It was defined so as to interpolate between the noncommutative
monomial basis introduced in \cite{Tev07} and the ribbon (Schur-like)
basis of noncommutative symmetric functions. Introducing
certain natural constraints led to a definition involving 
the special inversion statistic defined in \cite{NTW}. 
However, the relation between both approaches was unclear. 

\medskip
The aim of this paper is to 
clarify the relations between all these different approaches. 


 We shall uncover the relations between the contructions
of \cite{NTW} and of \cite{Tev} by first introducing a new multivariate
analogue of a classical automorphism of noncommutative symmetric functions, 
which will allow us to apply the methods of \cite{LNT} 
to find a recurrence for the matrices of \cite{Tev}, leading to a
natural multiparameter analogue, preserving a factorization property    
of the Kostka-like matrices.

The present paper started with the observation that
the Hall-Littewood basis of \cite{Tev} could
be obtained from one of the bases defined in \cite{NTW} by applying
a version of the so-called $(1-t)$-transform of noncommutative symmetric functions.
A second observation was that the matrices expressing these functions
in a suitably modified ribbon basis satisfied a recurrence relation of the same type
as those of \cite{LNT}. Since it has been shown in \cite{LNT} that such matrices
could be defined with many more parameters, we were naturally led to look
for a multiparameter version of the functions of \cite{Tev}. To this aim,
we had to find a multiparameter analogue  of the $(1-t)$-transform. It turns out
that such a map does exist. It is a morphism of algebras, though not
of coalgebras, which may be the reason for which it had been overlooked for a long time,
despite its simplicity. Its inverse, applied to a complete symmetric function,
yields the multiparameter Klyachko element 
of \cite{NCSF2}, which, as shown
by McNamara and Reutenauer \cite{mNR} does indeed reduce to a Lie idempotent
under an appropriate specialization (see also \cite{Dz,NT-col}).

The new recurrence for the (multivariate version of the) Kostka matrices of \cite{Tev}
allows us to provide a proof of a generalization of the product formula
annouced in this reference, as well as closed formulas for various other transition
matrices. All these results are easily proved by means of the Grassmann algebra
formalism of \cite{LNT}. The techniques of \cite{LNT} also allow us to introduce
a second family of parameters, so as to obtain Macdonald-like bases. For those,
we only describe some transition matrices. 

Our results do not provide multiparameter analogues of all the constructions
of \cite{NTW}, in particular, they do not seem to be related to the combinatorics
of permutation tableaux.
So, we conclude with an appendix sketching another approach to introducing
more parameters in the constructions of \cite{NTW}, by refining the special
inversion statistic with a code in the sense of \cite{HNT-Foata}.

{\small \bf Acknowlegements. -- \rm
Lenny Tevlin would like to thank his co-authors for their warm hospitality
at Marne-la-Vall\'ee, where this work was initiated.
}
\section{Notations}

Our notations for noncommutative symmetric functions will be as
in~\cite{NCSF1,NCSF2}. 
Here is a brief reminder.

The Hopf algebra of noncommutative symmetric
functions is denoted by $\Sym$, or by $\Sym(A)$ if we consider the realization
in terms of an auxiliary alphabet. Bases of $\Sym_n$ are labelled by
compositions $I$ of $n$. The noncommutative complete and elementary functions
are denoted by $S_n$ and $\Lambda_n$, and the notation $S^I$ means
$S_{i_1}\dots S_{i_r}$. The ribbon basis is denoted by $R_I$.
The notation $I\vDash n$ means that $I$ is a composition of $n$.
The conjugate composition is denoted by $I^\sim$. The length of $I$ is denoted by $\ell(I)$
and its weight by $|I|$. The mirror image of $I$ is denoted by $\bar I$.
If $I=(i_1,\ldots,i_r)$ is finer than $J=(j_1,\ldots,j_s)$, the
refining composition $I_J$ is the composition of $r=\ell(I)$ whose
$k$th part is the number of parts of $I$ composing $j_k$. For example,
if $I=(111122311)$ and $J=(3325)$, $I_J=(3213)$.

The product $R_IR_J=R_{I\triangleright J}+R_{IJ}$ of two ribbon Schur functions is the sum of
the two terms given by $I\triangleright J=(i_1,\ldots,i_r+j_1,\ldots,j_s)$
and $IJ=(i,_1,\ldots,i_r,j_1,\ldots,j_s)$.

The graded dual of $\Sym$ is $QSym$ (quasi-symmetric functions).
The dual basis of $(S^I)$ is $(M_I)$ (monomial), and that of $(R_I)$
is $(F_I)$. The \emph{descent set} of $I=(i_1,\dots,i_r)$ is
$\Des(I) = \{ i_1,\ i_1+i_2, \ldots , i_1+\dots+i_{r-1}\}$.

The monomial basis of \cite{Tev07} is denoted here by
$\Psi_I$ (instead of $M^I$), as in \cite{HNTT}. This basis
should not be confused with the basis $\Psi^I$ of \cite{NCSF1}.
 An important convention, followed in all the papers
of the NCSF series, is that an upper index denotes a multiplicative basis
built on a sequence of generators $Z_n$:
$Z^I=Z_{i_1}\cdots Z_{i_r}$. For a one-part composition $I=(n)$,
$Z^{(n)}=Z_n$, and the notation $Z_I$ is used only to denote a non-multiplicative
basis such that $Z_{(n)}=Z_n$.

\section{Noncommutative Hall-Littlewood functions}

This section provides background on some constructions which
will be simplified and generalized in the sequel.

\subsection{The Hall-Littlewood functions of \cite{Tev}}

The functions $P_I(t)$ of \cite{Tev} are, for $I=(n)$, 
\begin{equation}\label{defPn}
	P_n(t;A)=\frac{S_n((1-t)A)}{1-t}
\end{equation}
and for $I=(i_1,\ldots,i_r)$, given by the quasideterminant
\begin{equation}\label{defPI}
P_I(t;A) = \frac{(-1)^{r-1}}{[r]_t}
\begin{vmatrix}
P_{i_r} & [1]_t & 0 & \dots & 0& 0\\
P_{ i_{r -1} + i_r} & P_{i_{r - 1}} & [2]_t & \ldots & 0 & 0 \\
\vdots &  \vdots & \vdots & \vdots &  \vdots & \vdots \\
P_{i_2 + \ldots + i_r} & \ldots & \ldots & \ldots&  P_{i_2}& [r - 1]_t \\
\fbox{$P_{i_1 + \ldots + i_r}$} & \ldots & \ldots & \ldots &P_{i_1 +
i_2}& P_{i_1}
\end{vmatrix}
\end{equation}
where $[k]_t=1+t+\cdots+t^{k-1}$.

This
is not the definition given in \cite{Tev}, but this
 expression is equivalent to the recurrence
\begin{equation}
[r]_tP_I = P_{i_1}P_{i_2,\ldots,i_r}-P_{i_1+i_2}P_{i3,\ldots,i_r}+\cdots + (-1)^{r-1}P_{i_1+\cdots+i_r}\,.
\end{equation}
which is satisfied by the $P_I$ of \cite{Tev}.

The $Q$-basis is just \cite{Tev}
\begin{equation}
Q_I = (t;t)_r P_I\,.
\end{equation}

As already mentioned, the definition of the basis $P_I(t)$ in
\cite{Tev} is equivalent to (\ref{defPn}) and (\ref{defPI}).
This formulation is natural, as it amounts to replacing each
$\Psi_k$ in the quasideterminantal definition of the
monomial basis $\Psi_I$ by its natural $t$-analogue, and the integers
on the superdiagonal
by the corresponding $t$-integers. 
Then, by definition, $P_I(1)=\Psi_I$, and it is immediate that $P_I(0)=R_I$,
which is the exact analogue of the monomial/Schur specializations of the ordinary
Hall-Littlewood $P$-functions \cite{Mcd}. 

The main results of \cite{Tev} are (1) the product rule for $P_IP_J$
(Eq. (\ref{prodP}) below), and
(2) a combinatorial expression for  the analogues of the Kostka polynomials, expressing in the classical case
the $P$-expansion of Schur functions. If we define $K_{IJ}(t)$ by
\begin{equation}
R_J(A)=\sum_I K_{IJ}(t)P_I(t;A)
\end{equation}
then, the result of \cite{Tev} is
\begin{equation}
K_{IJ}(t)=\tilde D_I^J(t) = t^{\maj(I)}D_I^J(t^{-1})
\end{equation}
where $D_I^J(q)$ is the matrix defined in \cite[Prop. 3.7]{NTW} (see the Appendix of
the present paper).

\subsection{Comparison with the Hall-Littlewood functions of \cite{NTW}}

In \cite{NTW}, the parameter $q$ plays the role of $t^{-1}$ here. 
In this reference, a deformation $R_I(q)$ of the ribbon basis is introduced
by means of a nonassociative $q$-product on $\Sym$, itself
obtained by a linear projection of 
an associative $q$-product on the Hopf algebra$\WQSym$ (based on packed word,
or equivalently, set compositions, or surjections, see, {\it e.g.}, \cite{HNTT,NTW}).
The $R_J(q)$ are then expanded on a deformation $L_I(q)$ of the fundamental 
basis $L_I$ (see \cite{Tev07,HNTT}). The matrix $D_I^J(q)$
is defined by
\begin{equation}
R_J(q)=\sum_ID_I^J(q)\Psi_I(A)\,.
\end{equation}

Let us define
\begin{eqnarray}
\label{eq:fund}
\tilde R_J(t) &= \sum_I \tilde D_I^J(t)\Psi_I(A)\\
\tilde L_I(t) &= \sum_I k_{\bar I \bar J}(t)\Psi_I(A)\\
\tilde R_J(t) &=\sum_I \tilde F_I^J(t)\tilde L_I(t)
\end{eqnarray}
where $k_{IJ}(t)$ is the coefficient of $R_I(A)$ in Hivert's Hall-Littlewood
function $H_J(t;A)$, and  $F_I^J(t)$ is defined in \cite[Eq. (92)]{NTW}.

Then, the relations between the constructions of \cite{NTW} and \cite{Tev}
can be sumarized as follows:
\begin{proposition}
Let $\phi_t$ be the $\KK[t]$-linear endomorphism of $\Sym$ defined by $\phi_t(P_I(t))=\Psi_I$.
Then, $\phi_t(R_J)=\tilde R_J(t)$.
\end{proposition}

It is not clear whether all these families of  noncommutative
symmetric functions admit compatible multiparameter
analogues.
A multivariate analogue of $k_{IJ}(t)$ is defined in \cite{HLT}.
But to obtain multivariate versions of all the matrices of \cite{NTW},
we would need a multivariate analogue of $\tilde R_I(t)$
compatible with these multivariate $k_{IJ}$ (which is not
the case of $\RR(\t)$ defined below).
Numerical experiments indicate that it is unlikely that such
an analogue with good compatibility properties with the present multiparameter
$\PP$-functions could be defined.

However, interesting multivariate analogues of the matrices $D_I^J$ can be defined
by a different method, see
Section \ref{rubdrap}.

%


\section{A multivariate analogue of the $(1-t)$-transform}

The $(1-t)$-transform is an important automorphism of the algebra
of symmetric functions. On power sums, it is defined by
$p_n((1-t)X)=(1-t^n)p_n(X)$. This transformation appears for example in the
character formula for Hecke algebras, and it is an essential ingredient
of the theory of Hall-Littlewood functions \cite{Mcd},
 where a $t$-deformed scalar product is defined by
\begin{equation}
\<f(X),g(X)\>_t =\<f(X),g((1-t)X)\>\,.
\end{equation}
The classical Hall-Littlewood functions come in several flavors.
One first defines the functions $P_\lambda(t;X)$. Then, the functions
$Q_\lambda(t;X)$ (which are scalar multiples of the $P_\lambda$) are defined
as their adjoint basis for the $t$-deformed scalar product. The $Q'_\lambda(t;X)$
are defined as dual to the $P_\lambda$ for the undeformed scalar product (admitting
the Schur functions as an orthonormal basis). The relation between $Q$ and $Q'$ is
thus
\begin{equation}
Q_\lambda(t;X) = Q'_\lambda(t;(1-t)X)\,.
\end{equation}

Most of the known noncommutative analogues of Hall-Littlewood or Macdonald functions
admit multiparameter versions \cite{HLT,LNT}. It is known that
the $(1-t)$ transform can be defined in $\Sym$ \cite{NCSF2}. 
If we want to work
with multiparameter analogues of the  Hall-Littlewood functions, we need a multiparameter
version of the $(1-t)$-transform. This can be done as follows.

Recall from \cite{JNT} the expansion 
\begin{equation}\label{ttransR}
R_I\left(( 1-t)A)\right)=(-1)^{\ell(I)}\sum_{|J|=|I|, r=\ell(J)}(-1)^r
(1-t^{j_r})t^{\sum_{k\in{\mathcal A}(I,J)}j_k}S^J(A)
\end{equation}
where 
\begin{equation}\label{defAA}
\AA(I,J)=\{s<\ell(J)|j_1+\cdots +j_s\not\in\Des(I)\}.
\end{equation}
Let $\t=(t_i)_{i\ge 1}$.
We introduce the following multivariate version by ``lowering the exponents'':
\begin{equation}\label{defRR}
\RR_I(\t; A)=(-1)^{\ell(I)}\sum_{|J|=|I|, r=\ell(J)}(-1)^r
\left((1-t_{j_r})\prod_{k\in{\mathcal A}(I,J)}t_{j_k} \right)
S^J(A)
\end{equation}
For example,
\begin{eqnarray}
	\RR_3 &= (1-t_3)S^3 -(1-t_1)t_2S^{21}-(1-t_2)t_1S^{12}+(1-t_1)t_1^2S^{111}\,,\\
	\RR_{21} &=-(1-t_3)S^3+(1-t_1)S^{21}+(1-t_2)t_1S^{12}-(1-t_1)t_1S^{111}\,.
\end{eqnarray}

We can also define 
\begin{equation}
\SS^I(\t;A)=\sum_{J\le I}\RR_J(\t;A)\,.
\end{equation}

\begin{theorem}
The $\SS$-basis is multiplicative:
\begin{equation}
\SS^I(\t)\SS^J(\t)=\SS^{IJ}(\t)\,.
\end{equation}
Thus, $\RR_I$ is the image of $R_I$ by the automorphism
\begin{equation}
\theta_\t:\ S_n(A) \longmapsto \SS_n(\t;A)\,.
\end{equation}
\end{theorem}

\Proof Let $I\vDash m$ and $J\vDash n$.
It is sufficient to prove that $\RR_I\RR_J=\RR_{I\triangleright J}+\RR_{IJ}$.
Substituting the expressions given by (\ref{defRR}) in this product, we get
on the one hand
\begin{equation}
\RR_I\RR_J =  \sum_{I',J'}(-1)^{\ell(I)+\ell(J)-r-s}
\prod_{\substack{k\in\AA(I,I')}{l\in\AA(J,J')}}t_{i'_k}t_{j'_l}\cdot
(1-t_{i'r})(1-t_{j'_s}) S^{I'J'}
\end{equation}
(where $r=\ell(I')$ and $s=\ell(J')$),
and distributing the factor $(1-t_{i'_r})$, this can be rewritten as
\begin{equation}
\begin{split}
\sum_K & (-1)^{\ell(I)+\ell(J)-r}(1-t_{k_r})\prod_{p\in\AA(IJ,K)}t_{k_p}S^K\\
&+
\sum_K   (-1)^{\ell(I)+\ell(J)-r-1}(1-t_{k_r})\prod_{p\in\AA(I\triangleright J,K)}t_{k_p}S^K\,,
\end{split}
\end{equation}
where $K$ runs over compositions of $m+n$ such that $m\in\Des(K)$.
On the other hand, we see on (\ref{defRR})
that  $\RR_{I\triangleright J}+\RR_{IJ}$ is given by the same expression,
where this time $K$ runs over all compositions of $m+n$. But these extra terms
cancel, since if $m\not\in\Des(K)$, then $\AA(IJ,K)=\AA(I\triangleright J,K)$.
\qed

\begin{theorem}
The inverse of the automorphism $\theta_\t:\ S_n\mapsto \SS_n(\t)$ is
\begin{equation}
\theta_\t^{-1}:\ 
S_n\mapsto
\KK_n(\t;A) =\sum_{I\vDash n}\frac{\prod_{d\in \Des(I)}t_d}{(1-t_1)(1-t_2)\cdots (1-t_n)}R_I(A)\,,
\end{equation}
(the multiparameter Klyachko element already encountered in \cite{NCSF2,mNR,Dz,NT-col}).
\end{theorem}

\Proof
Let $((t))_n:=(1-t_1)\cdots (1-t_n)$.
 Substituting (\ref{defRR}) in the expression of $\KK_n$, we get
\begin{equation}
\frac1{((t))_n}\sum_{I,J\vDash n}(-1)^{\ell(I)+\ell(J)}
\prod_{d\in\Des(I)}t_d\prod_{k\in\AA(I,J)}t_{j_k}(1-t_{j_s})S^J
\end{equation}
so that the coefficient of $S^J$ is
\begin{equation}
\frac{(-1)^{s-1}}{((t))_n}\left(
\sum_{I\vDash n}\prod_{d\in\Des(I)}(-t_d)\prod_{k\in\AA(I,J)}t_{j_k}
\right) (1-t_{j_s})\,.
\end{equation}
For $J=(n)$, the sum in the parentheses is $((t))_n$, and for $s=\ell(J)>1$,
the sum vanishes since its terms cancel pairwise as follows.
If $p$ is a descent of $J$, and $D$ is a subset of $[n-1]$ not containing
$p$, there are exactly two compositions $I$ such that $\Des(I)\backslash \{p\}=D$,
and they have opposite coefficients in the sum.
\qed

Note that if we set $t_i=t$ for all $i$, then $\KK_n$ becomes the noncommutative
Eulerian polynomial $\AA_n^*(t;A)$ of \cite{NCSF1}.

\begin{example}{\rm
The entry $(I,J)$ in the following matrices is the coefficient
of $S^I$ in $\RR_J(\t)$:
\begin{equation} 
\begin{pmatrix}
1 - t_2       & t_2 - 1 \\
t_1 (t_1 - 1) & 1 - t_1
\end{pmatrix}
\qquad
\begin{pmatrix}
1 - t_3 & t_3 - 1 & t_3 - 1 & 1 - t_3 \\
t_2 (t_1 - 1) & 1 - t_1 & -t_2 (t_1 - 1) & t_1 - 1 \\
t_1 (t_2 - 1) & -t_1 (t_2 - 1) & 1 - t_2 & t_2 - 1 \\
-t_1^2 (t_1 - 1) & t_1 (t_1 - 1) & t_1 (t_1 - 1) & 1 - t_1
\end{pmatrix}
\end{equation}
The inverse matrices are
\begin{equation}
\begin{pmatrix}
\frac{1}{(1-t_1)(1-t_2)}   & \frac{1}{(1-t_1)^2} \\
\frac{t_1}{(1-t_1)(1-t_2)} & \frac{1}{(1-t_1)^2}
\end{pmatrix}
\qquad
\begin{pmatrix}
\frac{1}{(1-t_1)(1-t_2)(1-t_3)}    & \frac{1}{(1-t_1)^2(1-t_2)} &
\frac{1}{(1-t_1)^2(1-t_2)}  & \frac{1}{(1-t_1)^3} \\
\frac{t_2}{(1-t_1)(1-t_2)(1-t_3)}   & \frac{1}{(1-t_1)^2(1-t_2)} &
\frac{t_1}{(1-t_1)^2(1-t_2)} & \frac{1}{(1-t_1)^3} \\
\frac{t_1}{(1-t_1)(1-t_2)(1-t_3)}   & \frac{t_1}{(1-t_1)^2(1-t_2)} &
\frac{1}{(1-t_1)^2(1-t_2)}  & \frac{1}{(1-t_1)^3} \\
\frac{t_1t_2}{(1-t_1)(1-t_2)(1-t_3)} & \frac{t_1}{(1-t_1)^2(1-t_2)} &
\frac{t_1}{(1-t_1)^2(1-t_2)} & \frac{1}{(1-t_1)^3}
\end{pmatrix}
\end{equation}
}
\end{example}

\section{Some sequences of matrices}

We shall now introduce some simple sequences of matrices, which
will allow us to obtain directly a multiparameter version of the
noncommutative Hall-Littlewood functions of \cite{Tev}, as well
as their multiparameter Macdonald-like extension.

Let $\q=(q_n)_{n\ge 1}$ and $\t=(t_n)_{n\ge 1}$ be two sequences of commuting
indeterminates. For $n\ge 1$, we define three square matrices $A_n$, $B_n$, $T_n$
of size $2^{n-1}$, indexed by compositions of $n$ arranged in reverse lexicographic
order, {\it e.g.,} for $n=3$, by $3,21,12,111$ in this order.

\begin{definition}
For two compositions $I,J$ of $n$, let $\AA(I,J)$ be defined by
(\ref{defAA}). We define the matrix $A_n$ by
\begin{equation}
A_n(I,J)=\prod_{k\in\AA(\bar I,\bar J)}t_k\times\prod_{l\in\AA(\bar I^\sim,\bar J^\sim)}q_l\,.
\end{equation}
\end{definition}

The matrix $T_n=\diag(t_{\ell(I)})$ is diagonal. In particular, $T_1=(t_1)$.
We set $A_1=(1)$, and define $B_n$ as the result of substituting $q_{i+1}$
to $q_i$ in $A_n$. Thus, $B_1=(1)$. 

\begin{proposition}\label{recAn}
The matrices $A_n$ satisfy the recursion $A_1=(1)$ and for $n>1$,
\begin{equation}
	A_n = \begin{pmatrix}  
		B_{n-1} & A_{n-1}T_{n-1} \\
		q_1 B_{n-1} & A_{n-1}
	\end{pmatrix}
\end{equation}
\end{proposition}

\Proof 
Let us cut the matrix $A_n$ into four blocks. Each block is now indexed by pairs
of compositions of $n$ which can be represented by compositions of $n-1$.

For example, the top-left corner corresponds to compositions of the form
$(1\triangleright I, 1\triangleright J)$. In that case, we have
$ \AA(\overline{ 1\triangleright I},\overline{ 1\triangleright J}) =
\AA(\bar I,\bar J)$
and
$ \AA(\overline{ 1\triangleright I}^\sim,\overline{ 1\triangleright J}^\sim) =
\AA(\bar I^\sim,\bar J^\sim)+1$
since, for all compositions $K$,
$\overline{ 1\triangleright K}^\sim =  \bar K^\sim\cdot 1$.
So the matrix in the top-left corner of $A_n$ is obtained from $A_{n-1}$ by
substituting $q_{i+1}$ to $q_i$.

The same reasoning applies to the other blocks  of $A_n$.
\qed

For example, 
\begin{equation}
	A_2 = \begin{pmatrix}  
		1 & t_1 \\
		q_1  & 1
	\end{pmatrix}\,,
	\qquad\qquad
	A_3 = \begin{pmatrix}  
		1 & t_1 & t_1 & t_1t_2\\
		q_2  & 1 & q_1t_1 & t_2 \\
		q_1 & q_1t_1 & 1 & t_1 \\
		q_1q_2 & q_1 & q_1 & 1
	\end{pmatrix}
\end{equation}

\begin{equation}
A_4=
\begin{pmatrix}
1 & t_1 & t_1 & t_1 t_2 & t_1 & t_1 t_2 & t_1 t_2 & t_1 t_2 t_3 \\
q_3 & 1 & q_2 t_1 & t_2 & q_2 t_1 & t_2 & q_1 t_1 t_2 & t_2 t_3 \\
q_2 & q_2 t_1 & 1 & t_1 & q_1 t_1 & q_1 t_1 t_2 & t_2 & t_1 t_3 \\
q_2 q_3 & q_2 & q_2 & 1 & q_1 q_2 t_1 & q_1 t_2 & q_1 t_2 & t_3 \\
q_1 & q_1 t_1 & q_1 t_1 & q_1 t_1 t_2 & 1 & t_1 & t_1 & t_1 t_2 \\
q_1 q_3 & q_1 & q_1 q_2 t_1 & q_1 t_2 & q_2 & 1 & q_1 t_1 & t_2 \\
q_1 q_2 & q_1 q_2 t_1 & q_1 & q_1 t_1 & q_1 & q_1 t_1 & 1 & t_1 \\
q_1 q_2 q_3 & q_1 q_2 & q_1 q_2 & q_1 & q_1 q_2 & q_1 & q_1 & 1
\end{pmatrix}
\end{equation}

As for all known analogues of the $(q,t)$-Kostka matrices, we have: 
\begin{proposition}
The determinants
of the $A_n$ are products of linear factors:
\begin{equation}
\det A_n = \prod_{k=2}^n\prod_{i=1}^{n-1}(1-q_it_{k-i})^{{n-1\choose k-1}}\,.
\end{equation}
\end{proposition}

\Proof Subtracting $q_1$ times the upper half of the matrix to its lower half,
we have 
\begin{equation}
\begin{split}
	|A_n| = \left|\begin{matrix}  
		B_{n-1} & A_{n-1}T_{n-1} \\
		q_1 B_{n-1} & A_{n-1}
	\end{matrix}\right|
	&=
	 \left|\begin{matrix}  
		B_{n-1} & A_{n-1}T_{n-1} \\
		0 & A_{n-1}-q_1A_{n-1}T_{n-1}
	\end{matrix}\right|\\
&\\
	&=|B_{n-1}| |A_{n-1}| |I_{n-1}-q_1T_{n-1}|\,.
\end{split}
\end{equation}
\qed

For example,
\begin{equation}
|A_4|=
(1-q_1t_1)^3
(1-q_1t_2)^3
(1-q_2t_1)^3
(1-q_1t_3)
(1-q_2t_2)
(1-q_3t_1)\,.
\end{equation}

\section{Macdonald-like polynomials}

\subsection{An analogue of the $J$-basis}

We can now define a basis $\JJ_K(\q,\t;A)$ by interpreting the {\em columns} of $A_n$
as their expansion on the $\RR$-basis:
\begin{equation}
	\JJ_K(\q,\t;A) := \sum_{|I|=n}A_n(I,K)\RR_I(\t;A)\,.
\end{equation}
We regard it as an analogue of Macdonald's $J$-basis.

For example,
\begin{equation}
\begin{split}
\JJ_{31}(\q,\t) =&
t_1 \RR_4 + \RR_{31} + q_2 t_1 \RR_{22} + q_2 \RR_{211}\\& + q_1 t_1 \RR_{13}
+ q_1 \RR_{121} + q_1 q_2 t_1 \RR_{112} + q_1 q_2 \RR_{1111}.
\end{split}
\end{equation}

\subsection{Factorized expressions in the Grassmann algebra}

In \cite{LNT}, a general method for constructing multiparameter
bases of noncommutative symmetric functions with a Macdonald-like behaviour
has been described. The basic idea is to identify $\Sym_n$ with a Grassmann
algebra on $n-1$ variables $\eta_1,\ldots,\eta_{n-1}$, a ribbon $R_I$ being
encoded by the product $\eta_{d_1}\cdots\eta_{d_k}$, where $\Des(I)=\{d_1,\ldots,d_k\}$.

We need only a slight modification of the definitions of \cite{LNT}. Let
$U=(u_1,\ldots,u_{n-1})$ and $V=(v_1,\ldots,v_{n-1})$ be two sequences of
parameters. We set
\begin{equation}
\begin{split}
K_n(U,V)&=(u_1+v_1\eta_1)\cdots (u_{n-1}+v_{n-1}\eta_{n-1})\\
&=\sum_{I\vDash n}\prod_{d\in\Des(I)}v_d\prod_{e\not\in\Des(I)}u_e \ R_I\,.
\end{split}
\end{equation}

For each composition $I$ of $n$, build the pair of sequences 
$(U_I,V_I)=((u^I_j),(v^I_j))_{j=1}^{n-1}$ as follows from the ribbon diagram of $I$.
First, write $(1,q_1),\ldots,(1,q_k)$ in this order, starting from
the top left cell, in all cells which are non-descents of $I$.
Then, write $(t_1,1),\ldots,(t_l,1)$, in this order, in all cells
which are descents of $I$, starting from the bottom right cell, as
in the example below.

\small
\begin{equation}
(U_{4121},V_{4121})=
\moyyoung{(1,q_1) & (1,q_2)   & (1,q_3)   & (t_3,1) \cr
         \blank & \blank  & \blank  & (t_2,1) \cr
         \blank & \blank  & \blank  & (1,q_4) & (t_1,1)\cr
         \blank & \blank  & \blank  &\blank & \times\cr}
\end{equation}
\normalsize

\begin{theorem}
Let $\JJ'_I(\q,\t,A)=K_n(U_I,V_I)$. Then,
\begin{equation}
\JJ_I(\q,\t;A)=\theta_\t(\JJ'_I(\q,\t;A))\,.
\end{equation}
\end{theorem}

\Proof  By definition, the coefficient of $R_I$ in $\JJ'_J$ is
\begin{equation}
\prod_{d\in\Des(I)}v^J_d\prod_{e\not\in\Des(I)}u^J_e
=\prod_{k\in\AA(\bar I,\bar J)}t_k\times\prod_{l\in\AA(\bar I^\sim,\bar J^\sim)}q_l
=A_n(I,J)\,.
\end{equation}
\qed

As in \cite{LNT}, we also identify $QSym_n$ with a Grassmann algebra on dual
variables $\xi_i$, and set
\begin{equation}
\<\xi_D,\eta_E\>=\delta_{DE}\,.
\end{equation}
Then, if
\begin{equation}
L_n(X,Y) = (y_1-x_1\xi_1)\cdots (y_{n-1}-x_{n-1}\xi_{n-1})\,,
\end{equation}
we have
\begin{equation}
\<L_n(X,Y),K_n(U,V)\>=\prod_{i=0}^{n-1}(u_iy_i-v_ix_i)\,.
\end{equation}

\subsection{Specialization $\q=0$}
As announced,
under the specialization $\q=0$, the $\JJ$-functions  reduce to a multivariate version of
the Hall-Littlewood functions of \cite{Tev}, which as we shall see, are simply related
to those of \cite{NTW}.

Indeed, Definition 5 of \cite{Tev} reduces to (\ref{ttransR}) for $J=(n)$,
and expanding the quasi-determinant gives back the general case.

If we regard the $\JJ$-functions as analogues of the Macdonald $J$-functions,
we can  define natural analogues of the classical $P$ and $Q$-functions by
\begin{equation}
\prod_{i=1}^{\ell(I)}(1- t_i) \PP_I(\t;A)  = \QQ_I(\t;A) =\JJ_I(0,\t;A) 
\end{equation}
Note that $ \QQ_n(\t;A) = \RR_n(\t;A)$. 

\begin{theorem}
The $\PP$-functions satisfy the recurrence
\begin{equation}\label{recPP}
\frac{1-t_r}{1-t_1}   \PP_I = \PP_{i_1}\PP_{i_2,\ldots,i_r}-\PP_{i_1+i_2}\PP_{i_3,\ldots,i_r}+\cdots + (-1)^{r-1}\PP_{i_1+\cdots+i_r}\,.
\end{equation}
Equivalently, we have the quasideterminantal expression
\begin{equation}\label{qdetPPI}
\PP_I(\t;A) = (-1)^{r-1}\frac{1-t_1}{1-t_r}
\begin{vmatrix}
\PP_{i_r} & 1-t_1 & 0 & \dots & 0& 0\\
\PP_{ i_{r -1} + i_r} & P_{i_{r - 1}} & 1-t_2 & \ldots & 0 & 0 \\
\vdots &  \vdots & \vdots & \vdots &  \vdots & \vdots \\
\PP_{i_2 + \ldots + i_r} & \ldots & \ldots & \ldots&  \PP_{i_2}& 1-t_{r-1} \\
\fbox{$\PP_{i_1 + \ldots + i_r}$} & \ldots & \ldots & \ldots &\PP_{i_1 + i_2}& P_{i_1}
\end{vmatrix}
\end{equation}
\end{theorem}

\Proof 
By definition,
\begin{equation}\label{QQHH}
\QQ'_I:= \theta_\t^{-1}(\QQ_I) = K_n(U_I,0^{n-1})\,.
\end{equation}
For example, 
\begin{equation}
\QQ'_{4121}=(t_3+\eta_4)(t_2+\eta_5)(t_1+\eta_7)\,.
\end{equation}
Thus, the product of a one-part $\QQ'_i$ by an arbitrary $\QQ'_J$ is
\begin{equation}
\QQ'_i\QQ'_J=1\cdot (1+\eta_i)\cdot K^{[i]}_n(U_J,0^{n-1})
\end{equation}
where $ K^{[i]}_n$ is obtained from $K_n$ by replacing each
$\eta_k$ by $\eta_{k+i}$. Writing 
\begin{equation}
1+\eta_i=1-t_{\ell(J)}+t_{\ell(J)}+\eta_i\,,
\end{equation}
we obtain
\begin{equation}
\QQ'_i\QQ'_J= (1-t_{\ell(J)})\QQ'_{i\triangleright J} +\QQ'_{IJ}\,.
\end{equation}

This implies that the $\QQ'_I$ and hence
also the $\QQ_I$ satisfy the recursion
\begin{equation}\label{recQQ}
\begin{split}
\QQ_I = &\QQ_{i_1}\QQ_{i_2,\ldots,i_r}-(1-t_{r-1})\QQ_{i_1+i_2}\QQ_{i_3,\ldots,i_r}\\
&+(1-t_{r-2})(1-t_{r-1})\QQ_{i_1+i_2+i_3}\QQ_{i_4,\ldots,i_r}+\cdots
 + (-1)^{r-1}((t))_{r-1}\QQ_{i_1+\cdots+i_r}\,,
\end{split}
\end{equation}
which is equivalent to \eqref{recPP}. 
\qed

For example,
\begin{equation}
\begin{split}
\QQ'_{4121}=& 1\cdot (1+\eta_4)\cdot (t_2+\eta_7)\\
&- (1-t_3)\cdot 1\cdot (1+\eta_5)\cdot(t_1+\eta_7)\\
&+(1-t_2)(1-t_3)\cdot 1\cdot (1+\eta_7)\\
&+ (1-t_1)(1-t_2)(1-t_3)\cdot 1\\
&= (t_3+\eta_4)(t_2+\eta_5)(t_1+\eta_7)\,.
\end{split}
\end{equation}

Thus, we have:
\begin{corollary}
For $\t=(t,t^2,t^3,\ldots)$, 
\begin{equation}
\JJ_I(0,\t)=Q_I(t)\,.
\end{equation}
\end{corollary}
\qed

For example,
\begin{equation}
Q_{21}(t)=S^{21}((1-t)A)-(1-t)S^3((1-t)A) =\RR_{21}(t)+t\RR_3(t)\,.
\end{equation}

\subsection{$\JJ(\q,\t)$ on $S$}

A first remarkable property of the $\JJ$-basis is its $S$-expansion.
Indeed, 
as for the $Q$-functions of \cite{Tev},
the coefficients of the transition matrices are simple products
of linear factors and can be explicitly described. For $n=2,3$, these are
\begin{equation}
\left(
\begin{matrix}
 (1-t_2)(1-q_1) & -(1-t_1)(1-t_2) \\
-(1-t_1)(t_1-q_1) & (1-t_1)(1-t_1^2) \\
\end{matrix}
\right)
\end{equation}

\begin{equation}
{\footnotesize
\left(
\begin{matrix}
   (1-q_1)(1-q_2)(1-t_3) &-(1-q_1)(1-t_1)(1-t_3)
 &-(1-q_1)(1-t_1)(1-t_3) & (1-t_1)(1-t_2)(1-t_3) \\
  -(t_2-q_2)(1-q_1)(1-t_1) & (1-t_1t_2)(1-q_1)(1-t_1)
 & (t_2-q_1)(1-t_1)^2     & -(1-t_1t_2)(1-t_1)(1-t_2) \\
  -(t_1-q_1)(1-q_2)(1-t_2)  & (t_1-q_1)(1-t_1)(1-t_2)
 & (1-t_1^2)(1-q_1)(1-t_2) &-(1-t_1t_2)(1-t_1)(1-t_2) \\
   (t_1-q_1)(t_1-q_2)(1-t_1)  &-(t_1-q_1)(1-t_1^2)(1-t_1)
 &-(t_1-q_1)(1-t_1^2)(1-t_1) & (1-t_1t_2)(1-t_1^2)(1-t_1)
\end{matrix}
\right)
}
\end{equation}

Equivalently, these matrices describe the expansion of $\JJ'_J$ on the basis
$\theta_\t^{-1}(S^I)$. In the Grassmann formalism,
\begin{equation}
\theta_\t^{-1}(S^I)=\frac{1}{((t))_I}K_n(1^{n-1},T_I)
\end{equation}
where $T_I=(t^{(I)}_j)_{j=1}^{n-1}$ is such that 
$t^{(I)}_j=t_{j-\sum_{k<j}i_k}$ if
$j$ is not a descent of $I$, and $t_j^{(I)}=1$ otherwise.

Recall the following construction from \cite{LNT}.

We encode a composition $I$ with descent set $D$ by the boolean word
$u=(u_1,\dots,u_{n-1})$ such that $u_i=1$ if $i\in D$ and $u_i=0$
otherwise.

Let $ \y= \{ y_u \}$ be a family of indeterminates indexed by all boolean
words of length $\leq n-1$. For example, for $n=3$, we have the six parameters
$y_0, y_1, y_{00}, y_{01}, y_{10}, y_{11}$.

Let $u_{m\dots p}$ be the sequence $u_{m}u_{m+1}\dots u_{p}$
and define
\begin{equation}
 P_I := (1+ y_{u_1} \eta_1)
        (1+y_{u_{1\dots2}} \eta_2) \dots
        (1+y_{u_{1\dots n-1}} \eta_{n-1})
\end{equation}
or, equivalently,
\begin{equation}
\label{defP}
 P_I := K_n(1^{n-1}Y_I)
\qquad\text{with}\quad Y_I=( y_{u_1},y_{u_{1\dots2}},\dots,y_{u}) =:(y_k(I))\,.
\end{equation}

Similarly, let
\begin{equation}
Q_I := (y_{w_1} - \xi_1) (y_{w_{1\dots2}}- \xi_2) \dots
       (y_{w_{1\dots n-1}}- \xi_{n-1}) =: L_n(Y^I,1^{n-1})
\end{equation}
where 
$ w_{1\dots k} = u_1\dots u_{k-1}\,\overline{u_k}$ where
$\overline{u_k}=1-u_k$, so that
\begin{equation}
Y^I=( y_{w_1},y_{w_{1\dots2}},\dots,y_{w_{1\dots n-1}})=(y^k(I))\,.
\end{equation}

Then \cite[Prop. 4.1]{LNT}, the
 bases $(P_I)$ and $(Q_I)$ are dual to each other, up to normalization:
\begin{equation}
\<Q_I,P_J\>=\<L_n(Y^I,1^{n-1}),K_n(1^{n-1},Y_J)\>= \prod_{k=1}^{n-1}(y^k(I)-y_k(J))\,,
\end{equation}
which is indeed zero unless $I=J$.

Thus, in the notation of \cite[Eq. (31)]{LNT}
\begin{equation}
\theta_\t^{-1}(S^I)=\frac{1}{((t))_I}P_I,
\end{equation}
where the parameters $y_u$ of the binary tree are defined by
\begin{equation}
y_u=
\begin{cases}
1 &\text{if $u=u'1$},\\
t_j&\text{if $u=u'10^j$.}
\end{cases}
\end{equation}

Since $\JJ'_J=K_n(U_J,V_J)$, it is sufficient to expand a generic
$K_n(U,V)$ on this basis. We have clearly
\begin{equation}
\<L_n(Y^I,1^{n-1}),K_n(1^{n-1},Y_I)\>=\delta_{IJ}\prod_{k=1}^{n-1}(y^k(I)-y_k(I)).
\end{equation}

Thus, the dual basis of $\tilde P^I=\theta_\t^{-1}(S^I)$ is
\begin{equation}
\tilde Q_I = (1-t_{i_r})L_n(Y^I,1^{n-1}),
\end{equation}
so that the coefficient of $S^I$ in $\JJ_J$ is
\begin{equation}\label{eqJJ2SS}
(1-t_{i_r})\prod_{k=1}^{n-1}(u_k^J -y^k(I)v_k^J).
\end{equation}

\subsection{$\JJ$ on $S$ with more parameters}

Although  formula \eqref{eqJJ2SS} is completely explicit,
we can  better visualize the structure of these matrices
by introducing more parameters.
Let
\begin{equation}
\RR_I(\w,\x) := \sum_{J\vDash n} (-1)^{\ell(I)-r} (1-w_{r})  \prod_{k\in{\mathcal A}(I,J)} x_{j_k} S^J(A).
\end{equation}
(two sequences of parameters).

If we interpret the above matrices as describing a basis
$\JJ'(\q,\t)$ on $\RR(\w,\x)$,
then
$\JJ'(\q,\t)$ on $S$ still factorizes in the same manner:

\begin{equation}
\left(
\begin{matrix}
 (1-w_2)(1-q_1)   &-(1-t_1)(1-w_2) \\
-(1-w_1)(x_1-q_1) & (1-w_1)(1-t_1x_1) \\
\end{matrix}
\right)
\end{equation}

\begin{equation}
{\tiny
\left(
\begin{matrix}
   (1-q_1)(1-q_2)(1-w_3) &-(1-q_1)(1-t_1)(1-w_3)
 &-(1-q_1)(1-t_1)(1-w_3) & (1-t_1)(1-t_2)(1-t_3) \\
  -(x_2-q_2)(1-q_1)(1-w_1) & (1-t_1x_2)(1-q_1)(1-w_1)
 & (x_2-q_1)(1-t_1)(1-w_1) &-(1-t_1x_2)(1-t_2)(1-w_1) \\
  -(x_1-q_1)(1-q_2)(1-t_2)  & (x_1-q_1)(1-t_1)(1-w_2)
 & (1-t_1x_1)(1-q_1)(1-w_2) & -(1-t_2x_1)(1-t_1)(1-w_2) \\
   (x_1-q_1)(x_1-q_2)(1-w_1)  &-(x_1-q_1)(1-t_1x_1)(1-w_1)
 &-(x_1-q_1)(1-t_1x_1)(1-w_1) & (1-t_1x_2)(1-t_1x_1)(1-w_1)
\end{matrix}
\right)
}
\end{equation}



Let us now describe the coefficient $g_{JI}$ of $S^J$ in $\JJ_I$.

Let $K=I\wedge J$ be the composition whose descent set is $\Des(I)\cap\Des(J)$.
Let $\ell(I)=l$, $\ell(J)=m$ and $\ell(K)=n$.
Denote by $D$ and $D'$ the descent sets of $K_I$ and $K_J$.
Define
\begin{equation}
Z_1 = (1-w_{j_{m}}) \prod_{k=1}^{n-1} (1-t_{l-d_k} x_{J_{d'_k}})
       \prod_{s\in[1,l-1]\backslash D} (1-t_{l-s}).
\end{equation}
The conjugate mirror compositions
$I'=\overline{I}^\sim$ and $J'=\overline{J}^\sim$ encode
the complementary descent sets. 
Set $K'=I'\wedge J'$,
$D''=\Des(K'_{I'})$, and
\begin{equation}
Z_2 = \prod_{k=1}^{\ell(K')-1} (1-q_{d''_k})
\end{equation}
Define also
$E = [1,m-1]\backslash D'$ and $E'=[1,\ell(I')-1] \backslash D''$.
Observe that $E$ and  $E'$ have the same cardinality, which we denote by $e$.
Finally, set
\begin{equation}
Z_3 = \prod_{k=1}^{e} (x_{j_{e_k}} - q_{e'_k}).
\end{equation}

\begin{theorem}\label{JJ2S}
The coefficient $g_{JI}$ of $S^J$ in $\JJ_I$ is
\begin{equation}
g_{JI}= Z = (-1)^{\ell(I)-\ell(J)} Z_1 Z_2 Z_3\,.
\end{equation}
\end{theorem}
\qed

For example, with
$I = (13122)$ and $J=(221112)$,
we have $K=(4122)$,
so that $D=\{2,3,4\}$ and $D'=\{2,3,5\}$.

Thus, $Z_1 = (1-w_2) (1-t_3x_2)(1-t_2x_1)(1-t_1x_1) (1-t_4)$.

Next, $I'=(21321)$, $J'=(1251)$ and $K'=(351)$.
Hence, $D''=\{2,4\}$ and
$Z_2 = (1-q_2)(1-q_4)$.

Finally, $E=[1,4]$ et $E'=[1,3]$, thus
$Z_3 = (x_2-q_1)(x_1-q_3)$.

The final result is then
\begin{equation}
Z = - (1-w_2) (1-t_3x_2)(1-t_2x_1)(1-t_1x_1) (1-t_4) (1-q_2)(1-q_4)
(x_2-q_1)(x_1-q_3).
\end{equation}


\section{The Hall-Littlewood specialization}

Let us now continue the investigation of the specialization $\q=0$.

\subsection{Transition matrices $\SS$ to $\QQ$}

The matrices $M(\SS,\QQ)$ have a simple structure For example,

\begin{equation}
\left(
\begin{matrix}
1 & 1-t_1 \\
. & 1 \\
\end{matrix}
\right)
\end{equation}

\begin{equation}
\left(
\begin{matrix}
1 & 1-t_1 & 1-t_1 & (1-t_1)^2 \\
. & 1 & . & 1-t_2 \\
. & . & 1 & 1-t_1 \\
. & . & . & 1
\end{matrix}
\right)
\end{equation}

\begin{equation}
\left(
\begin{matrix}
1 & 1-t_1 & 1-t_1 & (1-t_1)^2  & 1-t_1 & (1-t_1)^2  & (1-t_1)^2 & (1-t_1)^3 \\
. & 1 & . & 1-t_2 & . & 1-t_2 & . & (1-t_2)^2 \\
. & . & 1 & 1-t_1 & . & . & 1-t_2 & (1-t_1)(1-t_2) \\
. & . & . & 1 & . & . & . & 1-t_3 \\
. & . & . & . & 1 & 1 - t_1 & 1-t_1 & (1-t_1)^2 \\
. & . & . & . & . & 1 & . & 1-t_2 \\
. & . & . & . & . & . & 1 & 1-t_1 \\
. & . & . & . & . & . & . & 1
\end{matrix}
\right)
\end{equation}

As usual, column $J$ gives $\SS^J(\t)$ as a linear combination of  $\QQ_I(\t)$.
These matrices can also be interpreted as giving the expansion 
of $S^J$ on the $\QQ'_I$.

The $S$-functions are given by
\begin{equation}
S^J=K_n(1^{n-1},X_J)\, \quad\text{where $x_i(J)=$} \begin{cases}
1 &\text{if $i\in\Des(J)$},\\
0&\text{otherwise}\,.
\end{cases}
\end{equation}

The dual basis of $\QQ'_I$ has a simple description:
\begin{theorem}\label{dualQ}
In terms of dual Grassmann variables, the dual basis of $\QQ'_I$ is
\begin{equation}
\GG_I = L_n(A_I,B_I) = (b_1(I)-a_1(I)\xi_1)\cdots (b_{n-1}(I)-a_{n-1}(I)\xi_{n-1})
\end{equation}
where
\begin{equation}
(a_j(I),b_j(I))=
\begin{cases}
(t_{c(j)},1) &\text{if $j\in\Des(I)$},\\
(-1,0)&\text{otherwise,}
\end{cases}
\end{equation}
where $c(j)=\ell(I)+1-r(j)$, $r(j)$ being the index of the row
containing the cell $j$ in the ribbon diagram. 
\end{theorem}

\Proof It is clear that $\<\GG_I,\QQ'_I\>=1$, since the only contribution
comes from the term of highest degree in $\QQ'_I$. It is also clear that
$\<\GG_I,\QQ'_J\>=0$ if $\Des(J)\not\supseteq\Des(I)$. In the remaining cases,
the greatest descent of $J$ which is not a descent of $I$ will contribute
a factor $t_k-t_k$, so that the result is indeed zero. \qed

The coefficient of  $\QQ_I(\t)$
in $\SS^J(\t)$ is  thus
\begin{equation}
\<\GG_{I}, S^{J}\>=\prod_k (b_k(I)x_k(J)-a_k(I))\,. 
\end{equation}

For example, the dual basis element corresponding to
\begin{equation}
\QQ'_{3122}=(t_3+\eta_3)(t_2+\eta_4)(t_1+\eta_6)
\end{equation}
is
\begin{equation}
\GG_{3122}=(1-t_4\xi_1)(1-t_4\xi_2)\xi_3\xi_4(1-t_2\xi_5)\xi_6(1-t_1\xi_7)
\end{equation}
so that the coefficient of $\QQ'=_{3122}$ in $S^{11312}=(1+\eta_1)(1+\eta_2)(1+\eta_5)(1+\eta_6)$
is $(1-t_4)^2(1-t_2)$.

\subsection{Multiplicative structure}

An interesting feature of most noncommutative analogues of the $P,Q$-functions
is that their product can be explicitly described. This is again the case here.

Let us recall  the product of the $P_I(t)$ is, as stated (along with a
few typos) in~\cite{Tev07}:

\begin{equation}\label{prodP}
\begin{split}
P_I(t) P_J(t)
= 
\sum_{K\le I} t^{\maj({I_K}^\sim)}
& \left( \tbin{\ell(K)+\ell(J)}{\ell(I)} P_{K\cdot J}(t)\right. \\
       &  +\left. \tbin{\ell(K)+\ell(J)-1}{\ell(I)} P_{K\triangleright J}(t) \right) \,,
\end{split}
\end{equation}
so that the product of the $Q_I$ is
\begin{equation}
\begin{split}
Q_I(t) Q_J(t)
= \sum_{K\leq I} t^{\maj({I_K}^\sim)}
 & \left(  \frac{[\ell(J)]_t!}{[\ell(J)+\ell(K)-\ell(I)]_t!} Q_{K\cdot J}(t) \right.\\
    & \left.   + \frac{[\ell(J)]_t!}{[\ell(J)+\ell(K)-\ell(I)-1]_t!} Q_{K\triangleright J}(t)\right)\,,
\end{split}
\end{equation}
with the convention that  $t$-factorials of negative numbers are zero.
Note that this amounts to removing the terms $Q_{C}$ such that $\ell(C)<\ell(I)$.

Now, in the case of multiple $t_i$, we have a more complicated but still elegant
formula.

Let us first define a sequence $S(I,K)$ associated with two compositions, $I$
finer than $K$.
Let $D=\{d_1<d_2<\ldots<d_p\}$ be the descent set of $I_K^\sim$.
Then, $S(I,K)=(d_p,d_{p-1},\ldots,d_1)$.

For example, with $I=(111122311)$ and $J=(3325)$, $I_K=(3213)$,
$I_K^\sim=(113211)$ and
$D=\{1,2,5,7,8\}$
so that $S(I,K)$ is $[8,7,5,2,1]$.
The $j$-th element of this sequence will be denoted by $s^{I,K}_j$.

Now, thanks to \eqref{QQHH},
it follows from Theorem \ref{dualQ} that
\begin{theorem}\label{prodQ}
\begin{equation}
\begin{split}
\QQ_I(\t) \QQ_J(\t)
= \sum_{K\leq I}
  \left( \prod_{j=1}^{\ell(I)-\ell(K)}
    \right.     &\left.  (t_{s^{I,K}_j} - t_{j+s^{I,K}_j+\ell(J)+\ell(K)-\ell(I)}\right)
          \QQ_{K\cdot J}(\t)\\
    &  \left.  + (1 - t_{\ell(J)})
         \prod_{j=1}^{\ell(I)-\ell(K)}
           (t_{s^{I,K}_j} - t_{j+s^{I,K}_j+\ell(J)+\ell(K)-\ell(I)-1})
          \QQ_{K\triangleright J}(\t) \right)
\end{split}
\end{equation}
with the convention that terms containing a $t_i$ with $i<0$ are zero.
Again, this amounts to removing the terms $\QQ_{C}$ such that $\ell(C)<\ell(I)$.
\end{theorem}
\qed

For example,
\begin{equation}
\begin{split}
\QQ_{211}(\t) \QQ_{21}(\t)
&= (t_1-t_3)(t_2-t_3) \QQ_{421}(\t) \\
&+ (t_2-t_4) \QQ_{3121}(\t) + (1-t_2)(t_2-t_3) \QQ_{331}(\t)  \\
&+ (t_1-t_3) \QQ_{2221}(\t) + (1-t_2)(t_1-t_2) \QQ_{241}(\t)  \\
&+ \QQ_{21121}(\t) + (1-t_2) \QQ_{2131}(\t).
\end{split}
\end{equation}

\begin{equation}
\begin{split}
\QQ_{111}(\t) \QQ_{111}(\t)
&= (t_1-t_4)(t_2-t_4) \QQ_{3111}(\t) + (1-t_3)(t_2-t_3)(t_1-t_3) \QQ_{411}(\t)  \\
&+ (t_2-t_5) \QQ_{21111}(\t) + (1-t_3)(t_2-t_4) \QQ_{2211}(\t)  \\
&+ (t_1-t_4) \QQ_{12111}(\t) + (1-t_3)(t_1-t_3) \QQ_{1311}(\t)  \\
&+ \QQ_{111111}(\t) + (1-t_3) \QQ_{11211}(\t).
\end{split}
\end{equation}

\begin{equation}
\begin{split}
\QQ_{1211}(\t) \QQ_{21}(\t)
&= (t_2-t_4)(t_3-t_4) \QQ_{4121}(\t) \\
&+ (1-t_3)(t_3-t_4) \QQ_{3221}(\t)  \\
&+ (t_1-t_3)(t_2-t_3) \QQ_{1421}(\t)  \\
&+ (t_3-t_5) \QQ_{31121}(\t) + (1-t_2)(t_3-t_4) \QQ_{3131}(\t)  \\
&+ (t_2-t_4) \QQ_{13121}(\t) + (1-t_2)(t_2-t_3) \QQ_{1331}(\t)  \\
&+ (t_1-t_3) \QQ_{12221}(\t) + (1-t_2)(t_1-t_2) \QQ_{1241}(\t)  \\
&+ \QQ_{111111}(\t) + (1-t_2) \QQ_{11211}(\t).
\end{split}
\end{equation}

\subsection{Limit cases}
Since $ P_I(t) $ of~\cite{Tev} interpolate between ribbon Schur (at $t=0$) and monomial 
(at $t=1$) bases, it is interesting to investigate these limits in the multivariate version.

Let ${\bf b} = (b_n)_{n \geq 1} $ be a sequence of commuting indeterminates. For $I=(i_1,\ldots,i_r)$,
define  a multivariate deformation of the monomial function $\Psi_I$ by
\begin{equation}
\Psi_I({\bf b};A) =(-1)^{r-1} \frac{b_1}{b_r}
\begin{vmatrix}
\Psi_{i_r} & 1 & 0 & \dots & 0& 0\\
\Psi_{ i_{r -1} + i_r} & \Psi_{i_{r - 1}} & \frac{b_2}{b_1}& \ldots & 0 & 0 \\
\vdots &  \vdots & \vdots & \vdots &  \vdots & \vdots \\
\Psi_{i_2 + \ldots + i_r} & \ldots & \ldots & \ldots&  \Psi_{i_2}& \frac{b_{r-1}}{b_1} \\
\fbox{$\Psi_{i_1 + \ldots + i_r}$} & \ldots & \ldots & \ldots & \Psi_{i_1 +
i_2}& \Psi_{i_1}
\end{vmatrix}
\end{equation}
and a deformation of the complete functions by
\begin{equation}
S_r({\bf b};A) = \frac{b_1}{b_r}
\begin{vmatrix}
\Psi_{1} &  - \frac{b_{r-1}}{b_1}  & 0 & \dots & 0& 0\\
\Psi_{2} & \Psi_{2} & \frac{b_{r -2}}{b_1}& \ldots & 0 & 0 \\
\vdots &  \vdots & \vdots & \vdots &  \vdots & \vdots \\
\Psi_{r-1} & \ldots & \ldots & \ldots&  \Psi_{1}& - 1\\
\fbox{$\Psi_{r}$} & \ldots & \ldots & \ldots & \Psi_{2}& \Psi_{1}
\end{vmatrix},
\end{equation}
where $ \Psi_n $ is a powers sum (of the first kind) as in~\cite{NCSF1}.

\begin{theorem}
Consider the specialization $ \t =(t^{b_1}, t^{b_2} , \ldots) $.
Then, for $ t \rightarrow 1 $, we have
\begin{align}
& \PP_I(1) = \Psi_I \left({\bf b} \right) \quad \text{ and } \\
& \PP_I(0) = R_I  \left({\bf b} \right), 
\end{align}
where 
$ R_I \left({\bf b} \right) = \sum_{J \preceq I} (-1)^{\ell(I) - \ell(J)} S^J \left({\bf b} \right) $ 
as usual. 
\end{theorem}

\Proof This follows from (\ref{qdetPPI}).\qed




\section{Appendix: another way to introduce multiple parameters in the functions of \cite{NTW}}\label{rubdrap}

The constructions of the present paper do not lead to nice
multiparameter analogues of all the matrices defined in \cite{NTW}.
We shall briefly describe here another approach to refining the constructions of \cite{NTW}.

The transition matrices of \cite{NTW} admit a combinatorial
description in terms of a statistic on words called \emph{special inversions}.

This statistic is defined on packed words, {\it i.e.}, words over the positive integers
whose support is an initial interval.
Let $u=u_1\cdots u_n$ be a packed word. We say that an inversion $u_i=b>u_j=a$
(where $i<j$ and $a<b$) is {\em special} if $u_j$ is the {\em rightmost}
occurence of $a$ in $u$. Let $\sinv (u)$ denote the number of special
inversions in $u$.
Note that if $u$ is a permutation, this coincides with its ordinary inversion
number.

The \emph{word composition} $\WC(w)$ of
a packed word $w$ is the composition
whose descent set is given by the positions of the last occurrences of each
letter in $w$.
Let $\DesC(w)$ be the usual descent composition of $w$.

For two compositions $I$ and $J$ of the same integer $n$,  let $W(I,J)$ be the set of packed
words $w$ such that
\begin{equation}
\WC(w)=I \qquad\text{and}\qquad \DesC(w)\fatter J.
\end{equation}
Then, the coefficients $C_I^J(q)$ of \cite{NTW} are
\begin{equation}
C_I^J(q) = \sum_{w\in W(I,J)} q^{\sinv(w)}
\end{equation}
and the $D_I^J(q)$ are
\begin{equation}\label{DIJq}
D_I^J(q) = \sum_{w\in W'(I,J)} q^{\sinv(w)}.
\end{equation}
where  $W'(I,J)$ is the set of packed words
$w$ such that
\begin{equation}
\WC(w)=I \qquad\text{and}\qquad \DesC(w)=J.
\end{equation}

The usual inversion number of a permutation can be refined into a 
list of integers (the Lehmer code) of which it is the sum. Something analogous can be done for
special inversions. As in \cite{HNT-Foata}, 
noncommutative generating functions for such codes can be given
in the form of flagged ribbon Schur functions.

If $\A=(A_{n_1}\supseteq A_{n_2}\supseteq\ldots\supseteq A_{n_r})$ is a flag of $r$ alphabets,
the flagged ribbon Schur function $R_I(\A)$, for a composition $I$ of length $r$,
is the sum of all semistandard fillings of the ribbon diagram of $I$ such that
only letters of $A_{n_i}$ appear in the $i$th row.

Let $A_i$ be the ordered alphabet $\{a_0,\dots,a_i\}$.

Then, \eqref{DIJq} (Formula (53) of \cite{NTW}) can be rewritten as
\begin{equation}
D_I^J(q) = \sum_{w\in W'(I,J)} \left.\left(\prod_{i\in\sinv(w)} a_i\right)\right|_{a_i=q^i}.
\end{equation}

Thus, 
\begin{equation}
D_I^J(q) = R_{J}(\A_{I,J})|_{a_i=q^i}
\end{equation}
where the flag of alphabets
 $\A_{I,J}$ is defined in the following way.
Let $I=(i_1,\dots,i_r)$ and $J=(j_1,\dots,j_s)$.
Then draw the ribbon diagram of $I$ and put dots into the cells corresponding
to the descents of $J$, and also into the last cell of $I$. 
Then $\A_{I,J}$ is the sequence of $s$ alphabets $A_{k_1},\dots,A_{k_s}$ where
$k_l$ is $r$ minus the row number of the $l$-th dotted cell of $I$.

For example, with $I=(3,1)$ and $J=(1,2,1)$, the sequence of alphabets is
$A_1,A_1,A_0$ since there are two dotted cells in the first row and one in the
second row:
\begin{equation}
\young{\cdot &   & \cdot \cr
         \blank & \blank  & \cdot \cr
         }
\end{equation}
Here are all the sequences of alphabets (where $(A_{k_1},\dots,A_{k_s})$ has
been replaced by $k_1\dots k_s$ to enhance readibility) corresponding to all
pairs of compositions of size $3$ and $4$.

\begin{equation}
\left(
\begin{matrix}
0 & .  & .  & . \\
0 & 10 & 10 & . \\
0 & .  & 10 & . \\
0 & 10 & 20 & 210
\end{matrix}
\right)
\end{equation}

\begin{equation}
\left(
\begin{matrix}
0 & .  & .  & .  & .  & .  & .  & . \\
0 & 10 & 10 & .  & 10 & 110& .  & . \\
0 & .  & 10 & .  & 10 & .  & .  & . \\
0 & 10 & 20 & 210& 20 & 210& 210& . \\
0 & .  & .  & .  & 10 & .  & .  & . \\
0 & 10 & 10 & .  & 20 & 210& 210& . \\
0 & .  & 10 & .  & 20 & .  & 210& . \\
0 & 10 & 20 & 210& 30 & 310& 320&3210
\end{matrix}
\right)
\end{equation}

Note that \emph{all} coefficients of the previous matrix are defined
as the flagged ribbons $R_J(\A_I)$, and that  some of them vanish: 
first, the ribbons such that
$\ell(J)>\ell(I)$ but also other ribbons, \emph{e.g.},
$R_{21}(A_{12,21}) = R_{21}(A_0,A_0)$ which has to be zero since the word
$000$ does not have a descent in position $2$.

Let now $I$ be a composition
and define
\begin{equation}
S^J(\A_{I,J}) := \sum_{J'\leq J} R_{J'}(\A_{I,J'}),
\end{equation}
then
the $S^J$ are multiplicative (as the usual ones).

This implies in particular that the equivalent of the
$C_I^J(q)$ (see~(29) of~\cite{NTW}) satisfy an induction similar to (50)
of~\cite{NTW}:
the whole matrix
$SP_n$ (see Section 3.5.1 of the same paper) has as analog the matrix of the
$S^J(\A_{I,J})$. This explains in a satisfactory way why all its
coefficients are products of binomial coefficients.
Indeed, if one sends $a_i$ to $q^i$,
\begin{equation}
S_n(1,q,\dots,q^{s}) = \qbin{s+n}{n}\,.
\end{equation}

One can also see the picture the other way round:
\begin{equation}
S_n(1,a_1,\dots,a_s)
\end{equation}
is the natural multivariate analog of the binomial $\binom{n+s}{n}=S^n(s+1)$
(in $\lambda$-ring notation).
Indeed, the induction on the $S_n$ is the natural analog of the classical
induction on binomials:
\begin{equation}
S_n(1,a_1,\dots,a_s)
= S_n(1,a_1,\dots,a_{s-1})
+ S_{n-1}(1,a_1,\dots,a_{s}) a_s.
\end{equation}

Setting $a_i =q_i$ (commuting parameters), we obtain multiparameter
versions of the Hall-Littlewood functions of \cite{NTW}. 
Their Macdonald-like extensions are not known.


\bigskip

\end{document}